\documentclass{amsart}

\usepackage[mathcal]{euscript}
\usepackage{amssymb, amsfonts, xypic}

\newtheorem{theorem}{Theorem}[section]
\newtheorem{lemma}[theorem]{Lemma}
\newtheorem{corollary}[theorem]{Corollary}
\newtheorem{proposition}[theorem]{Proposition}

\theoremstyle{definition}
\newtheorem{definition}[theorem]{Definition}
\newtheorem{example}[theorem]{Example}

\theoremstyle{remark}
\newtheorem{remark}[theorem]{Remark}

\numberwithin{equation}{section}

\newcommand{\Z}{\mathbb{Z}}
\newcommand{\bbA}{\mathbb{A}}

\newcommand{\C}{\mathcal{C}}
\newcommand{\A}{\mathcal{A}}
\newcommand{\D}{\mathcal{D}}

\newcommand{\HRR}{\textup{HRR}}
\newcommand{\Ar}{\textup{Ar}}
\newcommand{\Et}{\textup{Et}}
\newcommand{\id}{\textup{id}}
\newcommand{\pros}{\textup{pro}}
\newcommand{\pro}{\textup{pro-}}
\newcommand{\ind}{\textup{ind-}}
\newcommand{\op}{\textup{op}}

\DeclareMathOperator{\Hom}{Hom}
\DeclareMathOperator*{\colim}{colim}

\newcommand{\ol}{\overline}
\newcommand{\tl}{\tilde}

\newcommand{\map}{\rightarrow}

\newcommand{\dfn}{\textbf} 
\newcommand{\mdfn}[1]{\dfn{\mathversion{bold}#1}} 

\begin{document}

\title{Calculating Limits and Colimits in Pro-Categories}

\author{Daniel C. Isaksen}
\address{Department of Mathematics\\
University of Notre Dame\\
Notre Dame IN 46556, U.S.A}
\email{isaksen.1@nd.edu}
\thanks{The author was supported by an NSF Postdoctoral Fellowship.
The author thanks Dan Christensen for his constructive comments,
especially concerning Remark \ref{remark:other-limit} and
Section \ref{section:cocompact}.}

\date{June 6, 2001}

\begin{abstract}
We present some
constructions of limits and colimits in 
pro-categories.  These are critical tools in several applications.
In particular, certain technical arguments concerning strict pro-maps
are essential for a theorem about \'etale homotopy types.
Also, we show that cofiltered {\em limits} in pro-categories commute
with finite {\em colimits}.
\end{abstract}

\maketitle

\section{Introduction}

Pro-categories have found various uses over the past few decades, from
algebraic geometry \cite{SGA} \cite{AM} \cite{EF}
to shape theory \cite{CP} \cite{EH} \cite{MS} to
geometric topology \cite{CJS}
to applied mathematics \cite[Appendix]{CP}.
Unfortunately, the definition of a pro-category is subtle and
complicated.  The technical complexity of many papers using pro-categories
bears this out.
As an example of the subtlety of pro-categories, it was only recently
discovered \cite{prospace}
that the category $\pro \C$ is cocomplete whenever $\C$ is cocomplete.
(One must carefully interpret \cite[App. 4.3]{AM}, which discusses colimits
in pro-categories.  This result shows that $\pro \C$ is cocomplete when
$\C$ is small, which is useful for categories like pro-finite groups
but not useful for categories like pro-sets or pro-spaces.)

In several recent projects involving pro-categories \cite{duality}
\cite{prospace} \cite{ethtpy} \cite{proahss}, we have found it necessary
to compute various limits and colimits.  Since these 
different computations are similar, we have collected them together 
in this article.  We hope that these calculational tools will be useful
to others studying pro-categories.

We begin in Section \ref{section:prelim}
with a review of the necessary background on pro-categories.
We emphasize two important and well-known
facts about pro-objects because they are central to the techniques
in this paper.  First, every pro-object
is isomorphic to a pro-object indexed by a cofinite directed set.
Second, every map in a pro-category has
a level representation.

The category of arrows $\Ar(\pro \C)$ in a pro-category is equivalent
to the pro-category $\pro (\Ar \C)$ \cite{Meyer}.  
More generally, for any finite loopless category $A$,
the functor category
$(\pro \C)^{A}$ is equivalent to $\pro (\C^{A})$.
We introduce in Section \ref{section:level-replacements} a 
new reindexing result for a class of infinite diagrams.  
Unfortunately,
we cannot conclude that $(\pro \C)^{A}$ and $\pro (\C^{A})$ are equivalent
for our class of indexing diagrams $A$,
but we do have an essentially surjective functor
\[
\pro (\C^{A}) \map (\pro \C)^{A}.
\]
This means that 
every diagram in $\pro \C$ has a level
representation.
The functor is not an equivalence in general, but fortunately this is not
necessary for our purposes.
In Section \ref{subsctn:strict}, we specialize 
to a particular situation involving diagrams of strict
pro-maps.  This is critical for a theorem about 
hypercover descent for the \'etale topological type \cite{ethtpy}.

With this reindexing result in hand, we give an explicit description 
of cofiltered limits in pro-categories.
Since every limit can be rewritten in terms of finite limits of
loopless diagrams and cofiltered limits,
one can in principle 
describe an arbitrary limit if one can describe these two
special kinds.

Calculation of finite limits of loopless diagrams in pro-categories
is well-known \cite[App. 4.2]{AM}.  One finds a level representation
for the diagram 
and then takes the levelwise limit.
Calculation of cofiltered limits in pro-categories is more complicated.
The description of these limits in Theorem \ref{theorem:limit} 
is the essential result of this paper.  

This leads to the main applications.  
Given a class of maps $C$,
a pro-map is essentially of type $C$ if it has
a level representation by maps in $C$.
This kind of pro-map plays a significant role in many studies
involving pro-categories.  For example, it is central to the
abelian structure on the category $\pro \A$ for any abelian category $\A$
\cite{AM} \cite{Duskin}  \cite{essential} \cite{Verdier}.
Also, it is an important part of various closed model structures
on pro-categories \cite{duality} \cite{EH} \cite{prospace} \cite{proahss}.

We prove that the class of maps that are essentially of type $C$ is
closed under filtered limits for any class $C$.
This fact has several immediate 
corollaries.  The most interesting is that
filtered limits
are exact in the abelian category 
$\pro \A$ when $\A$ is any abelian
category. Thus, $\pro \A$ is very different from familiar
abelian categories such as abelian groups, $R$-modules, presheaves,
or sheaves.  The category $\pro \A$ has many of the properties of
the opposites of these familiar abelian categories, even though it
is not equivalent to any of them.

In the last section, we make some constructions of colimits in 
pro-categories.  We include these results because they are
similar to ideas in other parts of this paper and are likely to be
useful for computing concrete colimits in pro-categories.

We make one brief comment about ind-categories.  All the results of this
paper dualize because $\ind \C$ is equivalent to $(\pro \C^{\op})^{\op}$.
For example, we could describe filtered colimits in ind-categories.
We give no details because ind-categories occur less often in our
intended applications.
See \cite{duality} for one exception.

\section{Preliminaries on Pro-Categories} \label{section:prelim}

\begin{definition} 
\label{defn:pro}
For a category $\C$, the category \mdfn{$\pro \C$} has objects all 
cofiltering
diagrams in $\C$, and 
$$\Hom_{\pro \C}(X,Y) = \lim_s \colim_t \Hom_{\C}
     (X_t, Y_s).$$
Composition is defined in the natural way.
\end{definition}

A category $I$ is \dfn{cofiltering} if the following conditions hold:
it is non-empty and small;
for every pair of objects $s$ and $t$ in $I$,
there exists an object $u$ together with maps $u \map s$ and
$u \map t$; and for every pair of morphisms $f$ and $g$ with the
same source and target, there exists a morphism $h$ such that $fh$ equals
$gh$.
Recall that a category is \dfn{small} if it has only a set 
of objects and a set of morphisms.
A diagram is said to be \dfn{cofiltering} if its
indexing category is so.

Objects of $\pro \C$ are functors from
cofiltering categories to $\C$.  
We use both set theoretic and categorical
language to discuss indexing categories;
hence ``$t \geq s$'' and ``$t \map s$'' mean the same thing 
when the indexing category is actually a directed set.

The word ``pro-object'' refers to objects of pro-categories.
A \dfn{constant} pro-object is one indexed
by the category with one object and one (identity) map.
Let $c: \C \map \pro \C$ be the functor taking an object $X$ to the 
constant pro-object with value $X$.
Note that this functor makes $\C$ a full subcategory of
$\pro \C$.  The limit functor $\lim^{\C}: \pro \C \map \C$
is the right adjoint of $c$.  We shall always write this functor
as $\lim^{\C}$ to distinguish it from the functor
$\lim^{\pros}$, which is the limit internal
to the category $\pro \C$.

Let $Y: I \map \C$ and $X: J \map \C$ be arbitrary pro-objects.
We say that $X$ is 
\dfn{cofinal} in $Y$ if there is a cofinal functor
$F: J \map I$ such that
$X$ is equal to the composite $Y F$.
This means that for every $s$ in $I$,
the overcategory $F \downarrow s$ is cofiltered.
In the case when $F$ is an inclusion
of directed sets, $F$ is cofinal if and only if
for every $s$ in $I$ there
exists $t$ in $J$ such that $t \geq s$.
The importance of this definition is that $X$ is isomorphic to $Y$ in
$\pro \C$.

A \dfn{level representation} of a map
$f:X \map Y$ is:
a cofiltered index category $I$;
pro-objects $\tl{X}$ and $\tl{Y}$ indexed by $I$ and isomorphisms
$X \map \tl{X}$ and $Y \map \tl{Y}$;
and a collection of
maps $f_{s}:\tl{X}_s \map \tl{Y}_s$ for all $s$ in $I$
such that for all $t \map s$ in $I$, there is a commutative diagram
\[
\xymatrix{
\tl{X}_{t} \ar[d] \ar[r] & \tl{Y}_{t} \ar[d] \\
\tl{X}_{s} \ar[r]        & \tl{Y}_{s}       }
\]
and such that the maps $f_{s}$ represent a pro-map 
$\tl{f}: \tl{X} \map \tl{Y}$
belonging to a commutative square
\[
\xymatrix{
X \ar[d] \ar[r]^{f} & Y \ar[d] \\
\tl{X} \ar[r]_{\tl{f}} & \tl{Y}     }
\]
in $\pro \C$.
That is, a level representation is just a natural transformation 
such that the maps $f_{s}$ represent the element $f$ of
\[
\lim_{s} \colim_{t} \Hom_{\C} (X_{t}, Y_{s}) \cong
\lim_{s} \colim_{t} \Hom_{\C} (\tl{X}_{t}, \tl{Y}_{s}).
\]
Every map has a level representation 
\cite[Appendix 3.2]{AM} \cite{Meyer}.

More generally, suppose given any diagram $A \map \pro \C: a \mapsto X^{a}$.
A \dfn{level representation}
of $X$ is: a 
cofiltered index category $I$; a 
functor $\tl{X}: A \times I \map \C: (a, s) \mapsto \tl{X}^{a}_{s}$;
and isomorphisms $X^{a} \map \tl{X}^{a}$ such that 
for every map $\phi:a \map b$ in $A$, 
$\tl{X}^{\phi}$ is a level representation for $X^{\phi}$.
In other words, $\tl{X}$ is a uniform level representation for all the maps
in the diagram $X$.
Not every diagram of pro-objects has a level representation.

Suppose that $X: I \map \C$ and $Y: J \map \C$ are two pro-objects.
A \dfn{strict representation} \cite[p.~36]{EF} 
of a map
$f:X \map Y$ is
a functor $F: J \map I$
and a natural transformation
$\eta:X \circ F \map Y$ such that
the maps $\eta_{s}: X_{F(s)} \map Y_{s}$ represent the
element $f$ of
\[
\lim_{s} \colim_{t} \Hom_{\C} (X_{t}, Y_{s}).
\]

More generally, a \dfn{strict representation} of a diagram $X$ in
$\pro \C$ consists of
strict representations $(F^{\phi}, \eta^{\phi})$ for every map 
$\phi$ in $X$ such that
for every pair of composable maps $\phi$ and
$\psi$ in $X$,  the functor $F^{\psi \phi}$ equals $F^{\phi} F^{\psi}$
and the natural transformation $\eta^{\psi \phi}$ equals 
$\eta^{\psi} \circ \eta^{\phi} F^{\psi}$.

If $C$ is a class of objects, then a pro-object $X$ is 
\mdfn{of type $C$} if each $X_{s}$ belongs to $C$.
A pro-object is 
\dfn{essentially} of type $C$ if 
it is isomorphic to a pro-object of type $C$.

Similarly, a level representation $\tl{X} \map \tl{Y}$ of a map $X \map Y$
\mdfn{is of type $C$} (where
$C$ is a class of maps) if each $\tl{X}_{s} \map \tl{Y}_{s}$ belongs to $C$.
A map is \dfn{essentially} of type $C$
if it has a
level representation of type $C$.

\begin{definition} \label{definition:cofinite}
A category $A$ is \dfn{loopless} if it has no non-identity
endomorphisms.  
A category $A$ is \dfn{cofinite} if it is small, loopless, and 
for every object $a$ of $A$,
the set of arrows in $A$ with source $a$ is finite.
\end{definition}

A pro-object or level representation is \dfn{cofinite directed} 
if it is indexed by a cofinite directed set.

For every cofiltered category $J$, there exists a cofinite directed
set $I$ and a cofinal functor $I \map J$
\cite[Th. 2.1.6]{EH} 
(or \cite[Expos\'e\, 1, 8.1.6]{SGA}).  
Therefore, every
pro-object is isomorphic to a cofinite directed pro-object.
Similarly, every map has
a cofinite directed level representation.
Thus, it is possible to restrict the definition of a pro-object to
only consider cofinite directed sets as index categories.  
However, we find this unnatural.  Many general notions of pro-categories are
best expressed in terms of cofiltered categories, not in terms of 
cofinite directed sets.
For example, consider \cite{essential}, in which
the author assumes that all pro-objects are indexed by directed sets.
The construction of limits on pages 12--13 is technically correct, but it 
produces a pro-object that is not indexed by a directed set!

On the other hand, we find it much easier to work with
cofinite directed pro-objects in practice.  Thus, most of our results
start by assuming without loss of generality 
that a pro-object is indexed by a cofinite directed set.
Cofiniteness is critical because many arguments and constructions proceed
inductively.

\section{Level Representations and Cofiltered Limits} 
\label{section:level-replacements}

We study in this section
the question of replacing diagrams of pro-objects with 
level representations.  We then use level representations
to construct cofiltered limits.

One of the fundamental tools for studying pro-categories is the fact
that every morphism in a pro-category has a level 
representation.
In fact, every finite loopless diagram in a pro-category has 
a level representation
\cite[App. 3.3]{AM}.  The following result 
\cite{Meyer} is an elegant explanation of these level representation
principles.

\begin{theorem} \label{theorem:finite-level-replacement}
Let $A$ be any finite loopless category, and 
let $\C$ be any category.  Then the categories 
$\pro (\C^{A})$ and $(\pro \C)^{A}$ are equivalent.
\end{theorem}

When $A$ is the category with two objects and one non-identity morphism,
the theorem says that the 
categories $\Ar(\pro \C)$ and $\pro (\Ar \C)$
are equivalent.  Here $\Ar(\D)$ is the category whose objects are morphisms
in $\D$ and whose morphisms are commutative squares in $\D$.
This special case will be useful below in 
Corollaries \ref{corollary:map-limit} and \ref{corollary:map-retract}.

Beware that the pro-objects of \cite{Meyer} are indexed by cofiltered
categories that are not necessarily small.  On the other hand, all of
our indexing categories are small; this is an important distinction to
keep in mind.  Nevertheless, Theorem \ref{theorem:finite-level-replacement}
(and its consequences) are still true, as Remark
\ref{remark:not-faithful} indicates.

One important application of Theorem \ref{theorem:finite-level-replacement}
is the construction of finite limits in $\pro \C$ 
as stated in the next proposition from \cite[App. 4.2]{AM}.

\begin{proposition} \label{proposition:finite-limits}
Let $\C$ be a category with finite limits.
Let $A$ be a finite loopless category, and let
$X: A \map \pro \C$ be a diagram in $\pro \C$.  
Let $\tl{X}: A \times I \map \C$ be a level representation, where
$I$ is some cofiltered category, and
let $Z$ be the pro-object indexed by $I$
given by $Z_{s} = \lim^{\C}_{a} \tl{X}^{a}_{s}$.
Then $Z$ is isomorphic to $\lim^{\pros}_{a} X^{a}$.
\end{proposition}

Theorem \ref{theorem:finite-level-replacement}
(and therefore Proposition \ref{proposition:finite-limits}) also applies
when 
$\C$ is small and 
$A$ is an arbitrary finite category (possibly having 
non-trivial endomorphisms) \cite{Meyer}.
This is useful for 
categories such as pro-finite groups but not for pro-sets, pro-spaces,
pro-spectra, or pro-abelian groups.  Unfortunately, our applications involve
categories that are not small, so we 
do not proceed in this direction.
Rather, we generalize
in a different direction to the class of 
cofinite categories as described in Definition \ref{definition:cofinite}.

Every finite loopless category is cofinite, but other
infinite categories may also satisfy these conditions.  
The most important examples for our purposes are cofinite directed
sets.
We give more examples of cofinite categories
in Section \ref{section:colimits} when we study colimits of pro-objects.

\begin{theorem} \label{theorem:cofinite-level-replacement}
Let $A$ be a cofinite category, and 
let $X: A \map \pro \C: a \mapsto X^{a}$ 
be a diagram of pro-objects.
Then 
the diagram $X$ has a level representation.
\end{theorem}

\begin{proof}
We may assume that each pro-object $X^{a}$ is cofinite directed with
index set $I^{a}$.
Choose $I$ to be an arbitrary cofinite directed set with cardinality
greater than or equal to the cardinalities of every $I^{a}$; 
this will be the index set for the level representation $\tl{X}$.
Also choose 
arbitrary set surjections $h^{a}: I \map I^{a}$.

For every $a$, 
we build a new pro-object $\tl{X}^{a}$ by constructing a cofinal
function
$f^{a}: I \map I^{a}$ and letting $\tl{X}^{a}_{s}$ equal $X^{a}_{f^{a}(s)}$.

Fix an element $a$ of $A$.
Assume that the function $f^{b}$
has already been constructed on all indices $b$ for which there exists
a map $a \map b$ in $A$.  Then we may
proceed by induction because $A$ is cofinite.

We may define $f^{a}$ inductively since $I$ is cofinite.
Let $s$ be an index in $I$, and suppose that $f^{a}$ has already been
defined for $t < s$.
We choose $f^{a}(s)$ satisfying the following properties.  
This is possible because there are only finitely many conditions.

First, choose
$f^{a}(s)$ sufficiently large so that
$f^{a}(t) \leq f^{a}(s)$ for all $t <  s$.  This guarantees that $\tl{X}^{a}$
is a pro-object.  

Second, choose $f^{a}(s)$ large enough so that
$f^{a}(s) \geq h^{a} (s)$.  This guarantees
that $f^{a}$ is cofinal so that the natural map $X^{a} \map \tl{X}^{a}$ is 
an isomorphism.  

Third, choose $f^{a}(s)$ large enough so that for all 
maps $a \map b$ in $A$, there are 
maps $X^{a}_{f^{a}(s)} \map X^{b}_{f^{b}(s)}$ representing 
$X^{a} \map X^{b}$ such that the diagrams
\[
\xymatrix{
X^{a}_{f^{a}(s)} \ar[d] \ar[r] & X^{b}_{f^{b}(s)} \ar[d] \\
X^{a}_{f^{a}(t)} \ar[r]        & X^{b}_{f^{b}(t)}        }
\]
commute for all $t < s$.  This guarantees that $\tl{X}^{a} \map \tl{X}^{b}$
is a level representation.

Finally, choose $f^{a}(s)$ large enough so that for all pairs of arrows
$a \map b$ and $b \map c$ in $A$,
the diagram
\[
\xymatrix{
X^{a}_{f^{a}(s)} \ar[r] \ar[dr] & X^{b}_{f^{b}(s)} \ar[d] \\
& X^{c}_{f^{c}(s)}                   }
\]
commutes.
This guarantees that $\tl{X}$ is a level representation.

Note that the isomorphisms $X^{a} \map \tl{X}^{a}$ 
are given by representatives
\[
X^{a}_{f^{a}(s)} \map X^{a}_{f^{a}(s)} = \tl{X}^{a}_{s}
\]
that are identity maps.
It follows that the diagram
\[
\xymatrix{
X^{a} \ar[r] \ar[d] & X^{b} \ar[d] \\
\tl{X}^{a} \ar[r] & \tl{X}^{b}  }
\]
commutes for every map $a \map b$ in $A$ because
both compositions 
are given by representatives 
\[
X^{a}_{f^{a}(s)} \map X^{b}_{f^{b}(s)} = \tl{X}^{b}_{s}.
\]
\end{proof}

Objects of $\pro (\C^{A})$ are level representations of diagrams in
$\pro \C$.
This viewpoint gives us a functor
$F: \pro (\C^{A}) \map (\pro \C)^{A}$.  

\begin{corollary} \label{corollary:full-functor}
The functor $F$ is essentially surjective
in the sense that every object
of $(\pro \C)^{A}$ is isomorphic to an object in the image of $F$.
\end{corollary}

\begin{proof}
This follows immediately from 
Theorem~\ref{theorem:cofinite-level-replacement}.
\end{proof}

\begin{remark} \label{remark:not-faithful}
In general, the functor $F$ is not full and faithful unless $A$ is finite.
We briefly explain why.  Let $X$ and $Y$ be any two objects of
$\pro (\C^{A})$.  Then 
\[
\Hom (X, Y) =
   \lim_{s} \colim_{t} \int_{a} \Hom(X^{a}_{t}, Y^{a}_{s}),
\]
where the symbol $\int$ refers to the usual end construction for calculating
morphisms in a diagram category \cite[\S~IX.5]{MacLane}.
On the other hand,
\[
\Hom(FX, FY) =
  \int_{a} \lim_{s} \colim_{t} \Hom(X^{a}_{t}, Y^{a}_{s}).
\]
The limit involved in calculating the end
does not commute with the colimit unless it is a finite limit.
The end is a finite limit if and only if $A$ is a finite category.
\end{remark}

\subsection{Level Representations for Strict Diagrams}
\label{subsctn:strict}

The reindexing of Theorem \ref{theorem:cofinite-level-replacement}
is very general but also not canonical; one must make choices everywhere.
We next study a specific situation in which the reindexing is more natural.
This situation arises in \cite{ethtpy}.

Suppose given a strict representation of a diagram $X$ indexed by $A$.
Suppose that $A$ is cofinite and that the index category $I^{a}$
of each $X^{a}$ has finite limits.
If $I^{a}$ is a directed set,
then this means that every pair of elements in $I^{a}$ has a least upper
bound.

We now discuss a natural way of constructing the reindexing of
Theorem \ref{theorem:cofinite-level-replacement}.
Let $I$ be the product of the categories $I^{a}$ for which 
$a$ is the source of no non-identity maps.  The idea is that $I$ is 
the product of the indexing categories of all objects ``at the bottom''
of the diagram.  Note that $I$ is cofiltered since arbitrary products
preserve cofiltered categories.

For every $s = (s^{a})$ in $I$, define $\tl{X}^{a}_{s}$ to be
$X^{a}_{t}$, where 
\[
t = \lim_{\phi:a \map b} I^{\phi} (s^{b}).
\]
The limit is calculated in the category $I^{a}$.  The idea is that $t$
is the smallest common refinement of the elements $I^{\phi}(s^{b})$
in $I^{a}$.
It is straightforward to check that this definition satisfies all of 
the conditions of Theorem \ref{theorem:cofinite-level-replacement}.

\begin{example} \label{example:rigid-hypercover}
The main motivation for this discussion of strict representations
is that it is critical for the proof of the hypercover descent theorem for
\'etale toplogical types \cite{ethtpy}.

For every scheme $X$, let $\HRR(X)$ be the category of rigid hypercovers
of $X$ \cite[Prop.~4.3]{EF} \cite{ethtpy}.  This is the cofiltered
index category for the \'etale topological type $\Et X$ of $X$, which
is a pro-space.
Rigid limits \cite{ethtpy} give us finite limits
in $\HRR(X)$.  

For every scheme map $f:X \map Y$, rigid pullback \cite{ethtpy}
induces a functor $f^{*}:\HRR(Y) \map \HRR(X)$.  
For every rigid hypercover $U$ of $Y$, there is a canonical 
rigid hypercover map $f^{*} U \map U$.
This induces
a strict representation of the map $\Et X \map \Et Y$ of pro-spaces.

The rigid pullback functors are compatible under composition in the sense
that $(gf)^{*} = f^{*} g^{*}$.  Therefore,
every diagram of schemes induces a strict representation of a diagram
of pro-spaces.
Hence, we can apply the ideas of this discussion to reindexing
the \'etale topological types of cofinite diagrams of schemes.
\end{example}

\section{Cofiltered Limits} \label{section:limits}

We can now give a model for constructing cofiltered limits
in pro-categories.  
Start with a functor $X: B \map \pro \C$, where $B$ is a 
cofiltered index category.  
As mentioned in Section \ref{section:prelim},
there exists a cofinite directed set
$A$ and a cofinal functor $A \map B$.
We may replace $X$ with the composition $A \map B \map \pro \C$
because the limit does not change.
Since $A$ is a cofinite category, we can 
construct a level representation
$\tl{X}$ of $X$
as given in Theorem \ref{theorem:cofinite-level-replacement}.  
A product of cofiltered categories is again a cofiltered category, so
the functor $\tl{X}: A \times I \map \C$ can be viewed as a pro-object.

\begin{theorem} \label{theorem:limit}
The pro-object $\tl{X}$ is isomorphic to $\lim^{\pros}_{a} X^{a}$.
\end{theorem}

\begin{proof}
Because each $\tl{X}^{a}$ is isomorphic to $X^{a}$, it suffices to show
that $\tl{X}$ is isomorphic to $\lim^{\pros}_{a} \tl{X}^{a}$.  By direct
calculation with the definition of morphism sets for $\pro \C$,
$\tl{X}$ satisfies the required universal property.
\end{proof}

\begin{remark} \label{remark:other-limit}
The reindexing provided by 
Theorem \ref{theorem:cofinite-level-replacement} is not absolutely
necessary to construct cofiltered limits of pro-objects.
For example,
\cite{essential} contains another description of cofiltered limits.
This construction is suitable for proving some of the applications
below but not all of them.  One disadvantage of this alternative 
approach is that it produces pro-objects that are very complicated in
the sense that their index categories are far from being cofinite
directed sets.
On the other hand, the cofiltered limits constructed with
the help of Theorem \ref{theorem:cofinite-level-replacement} are
relatively small and computable.

For completeness, we repeat here the construction of \cite{essential}.
Let $X: A \map \pro \C: a \mapsto X^{a}$ 
be a diagram of pro-objects indexed by a
cofiltered category $A$ such that each $X^a$ is
indexed by a cofiltered category $I^a$.
Define a category $I$ as follows.  The objects of $I$ are pairs
$(s,a)$, where $s$ belongs to $I^{a}$.  That is,
the set of objects is the disjoint union of the sets $I^{a}$.
For each map $\phi : a \map b$ in $A$, we are given a map
$X^{\phi}: X^{a} \map X^{b}$, which is represented by a coherent
family of elements 
$X^\phi_{s}$ in $\colim_{t} \Hom_{\C }(X^{a}_{t},X^{b}_{s})$.
The set of maps from $(s,a)$ to $(t,b)$ in $I$ is 
defined to be the disjoint union over the maps $\phi : a \map b$ 
of the subsets of $\Hom _{\C }(X^{a}_{t},X^{b}_{s})$
consisting of those elements which equal $X^{\phi}_{s}$ in
$\colim_{t} \Hom_{\C }(X^{a}_{t},X^{b}_{s})$.

Composition in $I$ is defined by restricting the composition in $\C$.
Then $I$ is a cofiltered category, and 
there is a functor $I \map \C$ to $\C$ sending $(s,a)$
to $X^{a}_{s}$ and sending a map to itself.
It can be checked from the definitions that this pro-object is the 
limit of the given diagram.

The idea behind the above construction is that we take a diagram
whose objects are of the form $X^{a}_{s}$ and whose morphisms are
all possible maps $X^{a}_{s} \map X^{b}_{t}$ representing one
of the pro-maps $X^{a} \map X^{b}$ in the diagram $X$.
Note that even if $A$ and each $I^{a}$ are cofinite directed sets,
the category $I$ is in general not a directed set.
\end{remark}

\section{Essentially Levelwise Maps}

Our constructions leads to a result about the behavior of
cofiltered limits on pro-objects that are essentially of type $C$.

\begin{theorem} \label{theorem:object-limit}
Let $C$ be any class of objects of a category $\C$, and let 
$cC$ be the image of this class under the constant functor $c$.
The closure of $cC$ under isomorphisms and 
cofiltered limits is equal to
the class of pro-objects
in $\pro \C$ that are essentially of type $C$.
\end{theorem}

\begin{proof}
First we show that the class of pro-objects that are essentially of
type $C$ is closed under cofiltered limits.
Let $X$ be a cofiltered diagram of pro-objects, each of which is 
essentially of type $C$.  We may assume that each $X^{a}$ is actually
of type $C$.
We use the method of 
\cite[Th.~2.1.6]{EH}
to replace each $X^{a}$ by a cofinite directed pro-object that is
still of type $C$.  Then we use this same method to replace the
diagram $X$ by a cofinite directed diagram
with the same limit such that each
$X^{a}$ is still a cofinite directed pro-object of type $C$.

Now we use the construction of Theorem \ref{theorem:limit}.
Since each $X^{a}_{s}$ belongs to $C$, each $\tl{X}^{a}_{s}$ belongs to
$C$ 
(see the proof of Theorem \ref{theorem:cofinite-level-replacement}).
Since each $\tl{X}^{a}_{s}$ belongs to $C$, the pro-object $\tl{X}$
is of type $C$.
Therefore, $\lim^{\pros}_{a} X^{a}$ is essentially of type $C$.

Now we must show that every pro-object that is essentially of type $C$
is isomorphic to 
a cofiltered limit of objects in $cC$.  
It suffices to show that every pro-object of type $C$ is a cofiltered
limit of objects in $cC$.  Let $X$ be a pro-object such that each
$X_{s}$ belongs to $C$.  By direct calculation of morphism sets,
$X$ is isomorphic to $\lim^{\pros}_{s} cX_{s}$.  Also, each $cX_{s}$ belongs
to $cC$ by assumption.
\end{proof}

\begin{corollary} \label{corollary:map-limit}
Let $C$ be any class of maps in $\C$.  Then maps in $\pro \C$ 
essentially of type $C$ are
closed under cofiltered limits.
\end{corollary}

\begin{proof}
Apply Theorem~\ref{theorem:object-limit} to the category $\Ar \C$.
We use that the categories $\Ar(\pro \C)$ and $\pro (\Ar \C)$ are equivalent.
\end{proof}

This corollary gives several immediate results about particular 
classes of maps in pro-categories that are of interest.

\begin{corollary} \label{corollary:strict-structure}
Let $\C$ be a sufficiently nice model category \cite[\S 2.3]{EH}
or a proper model category.
The class of cofibrations in the strict model structure 
\cite[\S 3.3]{EH} \cite{proahss}
for $\pro \C$ is closed under cofiltered limits.
The class of fibrations in the strict model structure 
\cite[\S 3.8]{EH} for $\ind \C$ is closed under filtered colimits.
\end{corollary}

\begin{corollary} \label{corollary:homotopy-structure}
The classes of cofibrations in the $\pi_{*}$-model structure for
pro-spaces \cite{prospace}, in the $\pi_{*}$-model structure for
pro-spectra \cite{proahss}, and in the $\pi^{*}$-model structure
for pro-spectra \cite{duality} are all closed under cofiltered limits.
The class of fibrations in the $\pi_{*}$-model structure for
ind-spectra \cite{duality} is closed under filtered colimits.
\end{corollary}

The closure of cofibrations in the $\pi^{*}$-model structure for
pro-spectra under cofiltered limits is a technical necessity for
\cite{duality}. 

Having studied cofiltered limits of essentially levelwise maps,
we now proceed to retracts of such maps.

\begin{theorem} \label{theorem:object-retract}
Let $C$ be any class of objects of a category $\C$.  The class of
objects of $\pro \C$ that are essentially of type $C$ is
closed under retracts.
\end{theorem}

\begin{proof}
Let $Y$ be essentially of type $C$, 
and let $X$ be another pro-object with maps
$f:X \map Y$ and $g:Y \map X$ such that the composition
$gf$
is the identity. 
Consider the countable tower
\[
\xymatrix@1{
\cdots \ar[r] & X \ar[r]^{f} & Y \ar[r]^{g} & X \ar[r]^{f} & 
  Y \ar[r]^{g} & X. }
\]
Since $gf$ is the identity, the limit of this tower is isomorphic to $X$.
On the other hand, the limit is also isomorphic to the limit of the tower
\[
\xymatrix@1{
\cdots \ar[r] & Y \ar[r]^{fg} & Y \ar[r]^{fg} & Y.}
\]
By Theorem~\ref{theorem:object-limit}, this limit
is essentially of type $C$.
\end{proof}

We can now state the following corollary, which appears in
\cite[Prop.~12.1]{prospace}.  The proof there is philosophically the same
as the one here.

\begin{corollary} \label{corollary:map-retract}
Let $C$ be any class of maps of a category $\C$.  The class of 
maps of $\pro \C$ that are essentially of type $C$ is
closed under retracts.
\end{corollary}

\begin{proof}
Apply Theorem~\ref{theorem:object-retract} to the category
$\Ar \C$ and use that $\pro (\Ar \C)$ and $\Ar (\pro \C)$
are equivalent.
\end{proof}

\section{Finite Colimits and Cofiltered Limits}

It is a well-known and useful fact that filtered colimits commute
with finite limits in most familiar categories \cite[Th.~IX.2.1]{MacLane}.
Using our constructions of cofiltered limits given in 
Section \ref{section:limits},
we prove that the opposite is true for pro-categories.

\begin{theorem} \label{thm:limit-exact}
Cofiltered {\em limits} commute with finite {\em colimits} in
$\pro \C$.
\end{theorem}

\begin{proof}
Let $A$ be a cofiltered index category, and let $B$ be a finite
index category.  Suppose given a functor 
\[
X: A \times B \map \pro \C: (a, b) \mapsto X^{a,b}.
\]
In order to compute limits with respect to $A$, we may replace
$A$ with a cofinite cofiltered index category as described 
at the beginning of Section \ref{section:limits}.
In order to compute colimits with respect to $B$, we may replace
$B$ by a finite loopless category; this is the usual method for
rewriting any finite colimit in terms of finite coproducts and 
coequalizers.

Therefore, we may assume that both $A$ and $B$ are cofinite.  It
follows that $A \times B$ is cofinite, so we may assume that
$X$ is a level representation.  Thus, we have a functor
\[
X: A \times B \times I \map \C: (a, b, s) \mapsto X^{a,b}_{s}.
\]

Recall that finite colimits of pro-objects may be computed levelwise.
By Theorem \ref{theorem:limit} and direct computation,
both $\lim^{\pros}_{a} \colim^{\pros}_{b} X^{a,b}$ and 
$\colim^{\pros}_{b} \lim^{\pros}_{a} X^{a,b}$ are isomorphic to the
pro-object
\[
A \times I \map \C: s \mapsto \colim_{s} X^{a,b}_{s}.
\]
\end{proof}

\begin{remark}
Pro-categories are often equivalent to the opposites of more familiar
categories.  For example, the category of pro-finite abelian groups is
equivalent to the opposite of the category of torsion abelian groups
\cite[Prop. 7.5]{EF}, and the category of pro-finite $k$-vector spaces
is equivalent to the opposite of the category of $k$-vector spaces.
See \cite{duality} for an analogous stable homotopy theoretic statement.
More generally, 
the category $\pro \C$ is the opposite of the
pro-representable functors on $\C$ \cite[App. 4.2]{AM} 
\cite{Stauffer}.  When $\C$
is small (as in the case of finite abelian groups or finite
$k$-vector spaces), the pro-representable functors are just the
left exact functors.  

If the category $\pro \C$ were equivalent to the opposite of 
a more familiar category in which filtered colimits commuted with
finite limits, then
Theorem \ref{thm:limit-exact} would follow.
When $\C$ is not small, there is no clean
description of the pro-representable functors.
Even the categories of pro-abelian groups or pro-sets 
do not seem to be the opposites
of any particularly familiar categories.  
The isomorphism types
of cocompact pro-abelian groups or pro-sets do not form a set,
as shown in Theorem \ref{theorem:cocompact}.
However, most familiar categories
have only a set of isomorphism types of compact objects.  
\end{remark}

\section{Pro-Abelian Categories}

Let $\A$ be an abelian category.  Then $\pro \A$ is again an abelian
category \cite[App. 4.5]{AM} \cite{Duskin}.
The class of monomorphisms is equal to the class of
maps that have level representations by
monomorphisms.  The same is
true for the class of epimorphisms.

Now we have consequence of Theorem \ref{thm:limit-exact}.

\begin{theorem} \label{theorem:pro-lim-exact}
Let $\A$ be an abelian category.  Then cofiltered limits are exact in
the category $\pro \A$.  
\end{theorem}

\begin{proof}
Limits are always left exact.  Cofiltered limits are right exact
by Theorem~\ref{thm:limit-exact}.
\end{proof}

The surprising part of the previous theorem is that 
filtered limits are right exact.  This is
contrary to what occurs in most familiar abelian categories.  The
failure of epimorphisms to be closed under cofiltered limits
is the source of the higher $\lim^{i}$ functors.

Filtered colimits in 
the category $\pro \A$ need not be exact.
This is another
notable difference between pro-categories and most familiar
abelian categories.  
Thus we could define derived colimit functors $\colim^{i}$.
We expect that there exists a dual Mittag-Leffler condition
on filtered systems in $\pro \A$ that guarantees that the higher
derived colimits vanish.  We do not explore these ideas here since
we have no application for them.

\begin{example}
We give an example of the inexactness of filtered colimits in the
category pro-abelian groups.
Let $A$ be the free abelian group generated by the countable basis
$a_{0}, a_{1}, a_{2}, \ldots$.
For $n \geq m \geq 0$, let $A[m,n]$ be the subgroup of $A$
generated by $a_{m}, a_{m+1}, \ldots, a_{n}$.  For $m \geq 0$, let
$A[m, \infty)$ be the subgroup of $A$ generated by
$a_{m}, a_{m+1}, a_{m+2}, \ldots$.

Let $X$ be the pro-abelian group
\[
\cdots \map A[0, 2] \map A[0, 1] \map A[0,0].
\]
The structure maps of $X$ are the obvious projections.

Recall the constant functor $c: \C \map \pro \C$.
For $n \geq 0$, there is a map $f_{n}: cA[0,n] \map X$ of pro-abelian
groups.  This map is the element of 
\[
\lim_{m} \Hom(A[0,n], A[0,m]) \cong
\Hom(A[0,n], \lim_{m} A[0,m]) \cong
\Hom(A[0,n], \prod_{m} \Z a_{m})
\]
given by the obvious inclusion $A[0,n] \map \prod_{m} \Z a_{m}$.
Each map $f_{n}$ is a monomorphism in the category of pro-abelian groups
because it has a level representation 
\[
\xymatrix{
\cdots \ar[r] & A[0,n] \ar[r] \ar[d] & A[0,n] \ar[r] \ar[d] & 
  A[0,n] \ar[d] \\
\cdots \ar[r] & A[0, n+2] \ar[r] & A[0, n+1] \ar[r] & A[0,n]}
\]
in which each vertical map is a monomorphism.

Note that the maps $f_{n}$ are compatible with the inclusions
$cA[0,n] \map cA[0,n+1]$ in the sense that the diagram
\[
\xymatrix{
\cdots \ar[r] & cA[0, n-1] \ar[d]_{f_{n-1}} \ar[r] &
  cA[0, n] \ar[d]_{f_{n}} \ar[r] & 
  cA[0, n+1] \ar[d]_{f_{n+1}} \ar[r] & \cdots \\
\cdots \ar[r] & X \ar[r]_{=} & X \ar[r]_{=} & X \ar[r]_{=} & \cdots}
\]
of pro-abelian groups commutes.
The colimit of the first row is $cA$ because $A$ is the union of the
groups $A[0,n]$ and because $c$ commutes with colimits since it is 
a left adjoint.  The colimit of the second row is just $X$.
Thus, the colimit of the above diagram is a map $f:cA \map X$.
Now $f$ is a filtered colimit of the monomorphisms $f_{n}$,
but we claim that $f$ is not a monomorphism.

In order to show that $f$ is not a monomorphism,
we construct a non-zero map $g:Y \map cA$
such that the composition $fg$ is zero.
Let $Y$ be the pro-abelian group
\[
\cdots \map A[2, \infty) \map A[1, \infty) \map A[0, \infty).
\]
The structure maps of $Y$ are the obvious inclusions.
Let $g$ be the element of
\[
\colim_{n} \Hom(A[n,\infty), A)
\]
represented by any of the obvious inclusions $A[n,\infty) \map A$.
Then $g$ is not the zero map because 
each map $A[n,\infty) \map A$ is non-zero.
On the other hand, the composition $Y \map X$ is the zero element of
\[
\lim_{m} \colim_{n} \Hom(A[n,\infty), A[0,m])
\]
because for every $m$ and every $n > m$, 
the composition $A[n,\infty) \map A \map A[0,m]$
is zero.
\end{example}

\section{Cocompact Objects} \label{section:cocompact}

As another application of our explicit computation of cofiltered limits,
we describe completely in this section the cocompact objects of
$\pro \C$.  

Recall that an object $X$ of a category
is \dfn{compact} if for every filtered system $Y$, the map
\[
\colim_{a} \Hom (X, Y^{a}) \map 
  \Hom (X, \colim_{a} Y^{a})
\]
is an isomorphism.  We dualize this notion as follows.

\begin{definition} \label{definition:cocompact}
An object $X$ of a category is \dfn{cocompact} if for
every cofiltered system $Y$, the map
\[
\colim_{a} \Hom (Y^{a}, X) \map 
  \Hom ({\lim_{a}} Y^{a}, X)
\]
is an isomorphism.
\end{definition}

At first glance, this appears to be a strange definition.  Note that
$X$ is cocompact if and only if $X$ is compact in the opposite category.

\begin{theorem} \label{theorem:cocompact}
A pro-object is cocompact if and only if it is isomorphic to a 
constant pro-object.
\end{theorem}

\begin{proof}
First consider a constant pro-object $cX$.  Let $Y$ be an arbitrary
cofiltered system of pro-objects indexed by $A$.  
As described at the beginning of Section \ref{section:limits}, we may
assume that $A$ is cofinite.  Thus, we may take $Y$ to be a level
representation.
It follows by Theorem \ref{theorem:limit} and direct calculation that both
$\Hom_{\pros} (\lim^{\pros}_{a} Y^{a}, cX)$ and
$\colim_{a} \Hom_{\pros} (Y^{a}, cX)$ are equal to
$\colim_{a} \colim_{s} \Hom ( Y^{a}_{s}, X)$.
This finishes one direction of the theorem.

Now suppose that $X$ is a cocompact pro-object.  Note that $X$ is
isomorphic to $\lim^{\pros}_{s} cX_{s}$, so the map
\[
\phi: \colim_{a} \Hom_{\pros} (cX_{s}, X) \map 
  \Hom_{\pros} (X, X)
\]
is an isomorphism.  Let $f:cX_{s} \map X$ be any map such that
$\phi (f)$ is equal to the identity on $X$, and let 
$g: X \map cX_{s}$ be the natural map.

By definition of $\phi$, the composition $fg$ is the identity on $X$.
In order to show that $X$ and $cX_{s}$ are isomorphic, it suffices
to show that the composition $gf$ is the identity on $cX_{s}$.

Define a map
\[
\psi: \Hom_{\pros} (X, X) \map \colim_{a} \Hom_{\pros} (cX_{s}, X)
\]
by the formula $\psi (h) = fh$.  Because $fg$ is the identity, it
is easy to check that $\psi$ is the inverse of $\phi$.
Since $\psi (gf)$ and $\psi (\id)$ both equal $f$, it follows that
$gf$ is the identity map.
\end{proof}

\section{Colimits in Pro-Categories} \label{section:colimits}

We assume in this section that $\C$ is cocomplete.
The category $\pro \C$ is also cocomplete \cite[Prop.~11.1]{prospace},
but colimits in $\pro \C$ are difficult to compute in general.
See \cite{Blanc} for a special situation.
We study a different special
situation in which explicit computations are possible.
In particular, 
we compute colimits in $\pro\C$ indexed by cofinite categories.
At first glance, this does not seem to be especially useful, so we 
provide two examples.

\begin{example}
Consider the diagram for calculating the realization of a simplicial space
$X$.  This diagram has one object 
$X_{n} \otimes \Delta[n]$ for each $n \geq 0$ and one object
$X_{n} \otimes \Delta[m]$ for each simplicial operator 
$\phi: [m] \map [n]$.  The maps of the diagram are of 
two types.  The first type is of the form 
$\id \otimes \phi_{*}: X_{n} \otimes \Delta[m] \map X_{n} \otimes \Delta[n]$,
and the second type is of the form
$\phi^{*} \otimes \id: X_{n} \otimes \Delta[m] \map X_{m} \otimes \Delta[m]$.

The colimit of this diagram is the realization $|X|$ of $X$.  Note that
the diagram is cofinite because each object $X_{n} \otimes \Delta[n]$
is the source of zero non-identity maps 
and each object 
$X_{n} \otimes \Delta[m]$ is the source of two non-identity maps.
Therefore, the techniques of this section apply to 
calculating realizations of simplicial pro-objects.
\end{example}

\begin{example}
Consider a countable sequence 
\[
\xymatrix@1{
X_{0} \ar[r]^{f_{0}} & X_{1} \ar[r]^{f_{1}} & X_{2} \ar[r]^{f_{2}} & \cdots.}
\]
In order to calculate $\colim^{\pros} X$, 
we may take the colimit of the diagram
\[
\xymatrix{
X_{0} \ar[d]_{=} \ar[dr]_{f_{0}} & 
X_{1} \ar[d]_{=} \ar[dr]_{f_{1}} & 
X_{2} \ar[d]_{=} \ar[dr]_{f_{2}} & \cdots \\
X_{0} & X_{1} & X_{2} & \cdots,}
\]
which is cofinite.
Thus, the techniques of this section apply to calculating colimits 
of countable sequences.
\end{example}

Let $A$ be a cofinite category, 
and let $I$ be a cofinite directed set.
Recall that an arbitrary product of cofiltered categories is again
cofiltered.  Similarly, an arbitrary product of directed sets
is again a directed set, but infinite products 
do not preserve cofiniteness.
Consider the subset $K$ of $\prod_{A} I$ consisting of tuples
$(s_{a})$ such that $s_{a} \geq s_{b}$ when $a \geq b$.  Define 
a partial ordering on $K$ by $(s_{a}) \geq (t_{a})$ if $s_{a} \geq t_{a}$ for
all $a$.  This is the ordering that $K$ inherits as a subset
of $\prod_{A} I$.

\begin{lemma}
The set $K$ is directed.
\end{lemma}

\begin{proof}
Let $(s_{a})$ and $(t_{a})$ be any two elements of $K$.  We construct
a common refinement $(u_{a})$ by induction on $a$.
This is possible since $A$ is cofinite.
Suppose that $u_{b}$ has
already been determined for $b < a$ such that $u_{b} \geq s_{b}$ and
$u_{b} \geq t_{b}$.  Choose $u_{a}$ such that $u_{a} \geq s_{a}$,
$u_{a} \geq t_{a}$, and $u_{a} \geq u_{b}$ for all $b < a$.  This is 
possible because there are only finitely many conditions on $u_{a}$.
\end{proof}

Note that $K$ is not necessarily cofinite.

\begin{lemma} \label{lemma:cofinal1}
The inclusion $K \map \prod_{A} I$ is cofinal.
\end{lemma}

\begin{proof}
Given an element $(s_a)$ of $\prod_{A} I$, 
we must find an element $(t_{a})$ of $K$ such that $(t_{a}) \geq (s_{a})$.
Suppose that $t_{b}$ has already been chosen for $b < a$ such that
$t_{b} \geq s_{b}$.  Choose $t_{a}$ such that $t_{a} \geq s_{a}$ and
$t_{a} \geq t_{b}$ for all $b < a$.  This is possible since $A$ is cofinite.
\end{proof}

\begin{lemma} \label{lemma:cofinal2}
The forgetful functor $U^{a}: K \map I: (s_{a}) \mapsto s_{a}$ is cofinal.
\end{lemma}

\begin{proof}
The projection functor $\prod_A I \map I$ is cofinal.  By Lemma
\ref{lemma:cofinal1}, the functor $U^{a}$ is a composition of
two cofinal functors, so it is also cofinal.
\end{proof}

Let $X: A \map \pro \C$ be a
functor.  
By Theorem~\ref{theorem:cofinite-level-replacement},
we know that $X$ has a level representation $\tl{X}: A \times I \map \C$
for some cofinite directed set $I$.
We now define a functor $\ol{X}: A \times K \map \C$.  
Let $s = (s_{a})$ be an element of $K$.  For every object $a$ in $A$, 
define $\ol{X}^{a}_{s}$ to be $\tl{X}^{a}_{s_{a}}$.
For every map $a \map b$ in $A$, define 
$\ol{X}^{a}_{s} \map \ol{X}^{b}_{s}$ to be either composition
in the commuting square
\[
\xymatrix{
\tl{X}^{a}_{s_{a}} \ar[r] \ar[d] & \tl{X}^{b}_{s_{a}} \ar[d] \\
\tl{X}^{a}_{s_{b}} \ar[r] & \tl{X}^{b}_{s_{b}}.   }
\]
Note that the vertical maps are the structure maps of the pro-objects
$\tl{X}^{a}$ and $\tl{X}^{b}$.  
Here we use that $s_{a} \geq s_{b}$ since $s$ belongs
to $K$.  Also note that the horizontal maps come from the level representation
$\tl{X}^{a} \map \tl{X}^{b}$ of the map $X^{a} \map X^{b}$.

One can verify that $\ol{X}$ is indeed a functor by a straightforward 
diagram chase.

\begin{theorem} \label{theorem:pro-colimit}
The pro-object $Z$ indexed by $K$ given by 
$Z_{s} = \colim^{\C}_{a} \ol{X}^{a}_{s}$ is
the colimit in the category $\pro \C$ of the functor $X$.
\end{theorem}

\begin{proof}
Let $W$ be an arbitrary pro-object.
Then direct computation with the definition of morphism sets in
$\pro \C$ shows that
\[
\Hom_{\pros} (Z, W) \cong \lim_{a \in A} \Hom_{\pros} (X^{a}, W).
\]
This uses Lemma \ref{lemma:cofinal2} and 
the surprising fact that the canonical map
\[
\colim_{s \in K} \lim_{a \in A} \Hom_{\C} (\tl{X}^{a}_{s_{a}}, W_{t}) \map
\lim_{a \in A} \colim_{s \in K} \Hom_{\C} (\tl{X}^{a}_{s_{a}}, W_{t})
\]
is an isomorphism, which uses Lemma \ref{lemma:cofinal1}.
\end{proof}

\end{document}